\documentclass[12pt]{article}
\usepackage{amsfonts}
\usepackage{bbm}
\usepackage{mathrsfs}
\usepackage{amssymb}

\title{Graded Monomial Ordering for $\mathbb{N}$-graded and $\mathbb{N}$-filtered Solvable Polynomial Algebras of $({\cal B},d(~))$-type}

\author{Huishi Li\thanks{e-mail: huishipp@yahoo.com}\\
{\small Department of Applied Mathematics, College of Information Science and Technology}\\
{\small Hainan University,  Haikou 570228, China}}

\date{}

\begin{document}
\maketitle
\begin{center}
\begin{minipage}{135mm}
{\small {\bf Abstract.} Let $K$ be a field, and $A=K[a_1,\ldots
,a_n]$ a solvable polynomial algebra in the sense of [K-RW, {\it J. Symbolic
Comput.}, 9(1990), 1--26].  It is shown that if $A$ is an $\mathbb{N}$-graded  algebra of $({\cal B},d(~))$-type, the $A$ has a graded monomial ordering $\prec_{gr}$. It is also shown that $A$ is an $\mathbb{N}$-filtered algebra of $({\cal B},d(~))$-type if and only if $A$ has a graded momomial ordering, where ${\cal B}$ is the PBW basis of $A$.}
\end{minipage}\end{center} {\parindent=0pt\vskip 6pt

\vskip
6pt

}

\vskip .5truecm

\def\hang{\hangindent\parindent}
\def\textindent#1{\indent\llap{#1\enspace}\ignorespaces}
\def\re{\par\hang\textindent}

\def\v5{\vskip .5truecm}\def\QED{\hfill{$\Box$}}\def\hang{\hangindent\parindent}
\def\textindent#1{\indent\llap{#1\enspace}\ignorespaces}
\def\item{\par\hang\textindent}
\def \r{\rightarrow}\def\OV#1{\overline {#1}}
\def\normalbaselines{\baselineskip 24pt\lineskip 4pt\lineskiplimit 4pt}
\def\mapdown#1{\llap{$\vcenter {\hbox {$\scriptstyle #1$}}$}
                                \Bigg\downarrow}
\def\mapdownr#1{\Bigg\downarrow\rlap{$\vcenter{\hbox
                                    {$\scriptstyle #1$}}$}}
\def\mapright#1#2{\smash{\mathop{\longrightarrow}\limits^{#1}_{#2}}}
\def\NZ{\mathbb{N}}\def\mapleft#1#2{\smash{\mathop{\longleftarrow}\limits^{#1}_{#2}}}

\def\LH{{\bf LH}}\def\LM{{\bf LM}}\def\LT{{\bf
LT}}\def\KX{K\langle X\rangle} \def\KS{K\langle X\rangle}
\def\B{{\cal B}} \def\LC{{\bf LC}} \def\G{{\cal G}} \def\FRAC#1#2{\displaystyle{\frac{#1}{#2}}}
\def\SUM^#1_#2{\displaystyle{\sum^{#1}_{#2}}} \def\O{{\cal O}}  \def\J{{\bf J}}\def\BE{\B (e)}
\def\PRCVE{\prec_{\varepsilon\hbox{-}gr}}\def\BV{\B (\varepsilon )}\def\PRCEGR{\prec_{e\hbox{\rm -}gr}}

\def\KS{K\langle X\rangle}
\def\LR{\langle X\rangle}\def\T#1{\widetilde #1}
\def\HL{{\rm LH}}\def\NB{\mathbb{N}}\def\HY{\hbox{\hskip .03truecm -\hskip .03truecm}}\def\F{{\cal F}}

\def\hang{\hangindent\parindent}
\def\textindent#1{\indent\llap{#1\enspace}\ignorespaces}
\def\re{\par\hang\textindent}

\def\hang{\hangindent\parindent}
\def\textindent#1{\indent\llap{#1\enspace}\ignorespaces}
\def\item{\par\hang\textindent}
\def\normalbaselines{\baselineskip 24pt\lineskip 4pt\lineskiplimit 4pt}

\def\v5{\vskip .5truecm}
\def\QED{\hfill{$\Box$}}

\def\mapdown#1{\llap{$\vcenter {\hbox {$\scriptstyle #1$}}$}
                                \Bigg\downarrow}
\def\mapdownr#1{\Bigg\downarrow\rlap{$\vcenter{\hbox
                                    {$\scriptstyle #1$}}$}}

\def\mapright#1#2{\smash{\mathop{\longrightarrow}\limits^{#1}_{#2}}}
\def\mapleft#1#2{\smash{\mathop{\longleftarrow}\limits^{#1}_{#2}}}
\def \r{\rightarrow}

\def\FRAC#1#2{\displaystyle{\frac{#1}{#2}}}
\def\SUM^#1_#2{\displaystyle{\sum^{#1}_{#2}}}

\def\OV#1{\overline {#1}}
\def\T#1{\widetilde {#1}}

\def\NZ{\mathbb{N}}
\def\MB{\mathbb{B}}

\def\LH{{\bf LH}}
\def\LM{{\bf LM}}
\def\LT{{\bf LT}}
\def\LC{{\bf LC}}

\def\B{{\cal B}}
\def\O{{\cal O}}
\def\G{{\cal G}}
\def\F{{\cal F}}
\def\N{{\cal N}}
\def\S{{\cal S}}

\def\VF{\varphi}
\def\BE{\B (e)}
\def\VE{\varepsilon}
\def\PRCVE{\prec_{\varepsilon\textrm{\tiny -}gr}}
\def\BV{\B (\varepsilon )}

\def\PRCEGR{\prec_{e\textrm{\tiny -}gr}}
\def\PRECSVE{\prec_{s\textrm{\tiny -}\varepsilon}}
\def\PRECEVE{\prec_{e\textrm{\tiny -}\VE}}

\def\KX{K\langle X\rangle}
\def\GR{Gr\"obner~}


This note is a complement to [Li6] in which the graded monomial ordering for $\mathbb{N}$-graded and $\mathbb{N}$-filtered solvable polynomial algebras of $({\cal B},d(~))$-type is necessarily employed but its  existence is not established, and the structures of $\mathbb{N}$-graded and $\mathbb{N}$-filtered solvable polynomial algebras of $({\cal B},d(~))$-type are specified but not formally  defined.

\section*{1. $\NZ$-graded and $\NZ$-filtered Algebras of $(\B ,d(~))$-type}

Let $A=K[a_1,\ldots ,a_n]$ be a finitely generated $K$-algebra with the set of generators $\{ a_1,\ldots ,a_n\}$ and the PBW basis $\mathcal{B} =\{ a^{\alpha}=a_1^{\alpha_1}\cdots
a_n^{\alpha_n}~|~\alpha =(\alpha_1,\ldots ,\alpha_n)\in\NZ^n\}$. We first explore the sufficient and necessary condition that makes $A$ into  an $\NZ$-graded algebra such that \par
(a) the degree-0 homogeneous part of $A$ is equal to $K$;\par

(b) every $a^{\alpha}\in\B$ is a homogeneous element of $A$.\v5

Suppose that $A$ is an $\NZ$-graded algebra. Then by the definition of a graded algebra we know that $A=\oplus_{p\in\NZ}A_p$ with the degree-$p$
homogeneous part $A_p$ which is a subspace of $A$, such that $A_{p}A_{q}\subseteq A_{p+q}$ for all $p,q\in\NZ$. Furthermore, we assume that $A$ has the above  properties (a) and (b). Then $a_i\in A_{m_i}$ for some $m_i>0$, $1\le i\le n$. Thus,
writing $d_{\rm gr}(h)$ for the {\it
graded-degree} (abbreviated to {\it gr-degree}) of a nonzero\index{graded degree}\index{gr-degree}
homogeneous element $h\in A_p$, i.e., $d_{\rm gr}(h)=p$, we have $d_{\rm gr}(a_i)=m_i$, $1\le i\le n$.
It turns out that if $a^{\alpha}=a_1^{\alpha_1}\cdots a_n^{\alpha_n}\in\B$, then $d_{\rm
gr}(a^{\alpha})=\sum^n_{i=1}\alpha_im_i$. This shows that  $d_{\rm
gr}(~)$ gives rise to a {\it positive-degree function}\index{positive-degree function} on $\mathcal{B}$ (or equivalently on $A$) with respect to the $n$-tuple $(m_1,\ldots ,m_n)\in\NZ^n$, such that{\parindent=1.32truecm\par

\item{(1)} for every $p\in\NZ$, $A_p=K\hbox{-span}\{ a^{\alpha}\in\mathcal{B}~|~d_{\rm
gr}(a^{\alpha})=p\}$ (note that this is a finite dimensional subspace of $A$);\par

\item{(2)} for $1\le i<j\le n$, if $a_ja_i\ne 0$ and $a_ja_i=\sum^t_{i=1}\lambda_ia^{\alpha (i)}$, then
$d_{\rm gr}(a^{\alpha (i)})=m_j+m_i$, where $\lambda_i\in K^*$ and $a^{\alpha (i)}\in\mathcal{B}$ .\par}

Conversely, given an $n$-tuple $(m_1,\ldots ,m_n)\in\NZ^n$ with all $m_i>0$, there is a positive-degree function $d(~)$ on $A$ such that
$d(a^{\alpha})=\sum_{i=1}^n\alpha_im_i$ for all $a^{\alpha}=a_1^{\alpha_1}\cdots a_n^{\alpha_n}\in\B$, in particular $d(a_i)=m_i$, $1\le i\le n$, and for a nonzero $f=\sum^t_{s=1}\lambda_sa^{\alpha (s)}$ with $\lambda_s\in K^*$ and $a^{\alpha (s)}\in\B$, $d(f)=\max\{ d(a^{\alpha (s)})~|~1\le s\le t\}$.
With respect to $d(~)$, the vector space $A$ is equipped with an $\NZ$-graded structure $A=\oplus_{p\in\NZ}A_p$, where for every $p\in\NZ$,
$$A_p=K\hbox{-span}\{ a^{\alpha}\in\mathcal{B}|~d(a^{\alpha})=p\}$$ which is a finite dimensional subspace of $A$, in
particular, $A_0=K$. Clearly, every $a^{\alpha}\in\mathcal{B}$ is a homogeneous element of $A$.  However, in general this $\NZ$-gradation does not necessarily satisfy the condition
$$A_pA_q\subseteq A_{p+q},\quad p, q\in\NZ ,\leqno{(2')}$$ namely $A$ is not necessarily an $\NZ$-graded algebra with respect to the  $\NZ$-gradation determined by the PBW basis $\B$ of $A$ and the given positive-degree function $d(~)$ on $A$ (see Example given after Definition 1.4). To remedy this problem, considering the property (2) presented above, it is straightforward to verify that if
 $$\begin{array}{l} a_ja_i\ne 0,~a_ja_i=\sum^t_{s=1}\lambda_sa^{\alpha (s)}~\hbox{implies}~
 d(a^{\alpha (s)})=m_j+m_i, \\
 \hbox{where}~ \lambda_s\in K^*,~a^{\alpha (s)}\in\mathcal{B},~1\le i<j\le n,
\end{array}$$
then the condition ($2'$) is satisfied and consequently, $A$ is turned into an $\NZ$-graded algebra.\par

Summing up, the above discussion has led to  the following result and defination.{\parindent=0pt\v5

{\bf 1.1. Proposition}  Let $A=K[a_1,\ldots ,a_n]$ be a finitely generated $K$-algebra with the set of generators $\{ a_1,\ldots ,a_n\}$ and the PBW basis $\mathcal{B} =\{ a^{\alpha}=a_1^{\alpha_1}\cdots
a_n^{\alpha_n}~|~\alpha =(\alpha_1,\ldots ,\alpha_n)\in\NZ^n\}$. The following two statements are equivalent.}\par

(i) $A=\oplus_{p\in\NZ}A_p$ is an $\NZ$-graded algebra with the degree-$p$
homogeneous part $A_p$ which is a subspace of $A$, such that {\parindent=2truecm\par

\item{(a)} $A_0=K$;\par

\item{(b)} if $a^{\alpha}\in\mathcal{B}$, then $a^{\alpha}\in A_q$ for some $q\in\NZ$.\par}

(ii) With respect to a certain $n$-tuple $(m_1,\ldots ,m_n)\in\NZ^n$ in which all $m_i>0$, there is a positive-degree function $d(~)$ on $A$ such that
$d(a^{\alpha})=\sum_{i=1}^n\alpha_im_i$ for all $a^{\alpha}=a_1^{\alpha_1}\cdots a_n^{\alpha_n}\in\B$, in particular $d(a_i)=m_i$, $1\le i\le n$, and such that for $1\le i<j\le n$, if $a_ja_i\ne 0$ and $a_ja_i=\sum^t_{s=1}\lambda_sa^{\alpha (s)}$, then
$d(a^{\alpha (s)})=m_j+m_i$, where $\lambda_i\in K^*$ and $a^{\alpha (s)}\in\mathcal{B}$.\par\QED{\parindent=0pt\v5

{\bf 1.2. Definition}  Let $A=K[a_1,\ldots ,a_n]$ be a finitely generated $K$-algebra with the set of generators $\{ a_1,\ldots ,a_n\}$ and the PBW basis $\mathcal{B} =\{ a^{\alpha}=a_1^{\alpha_1}\cdots
a_n^{\alpha_n}~|~\alpha =(\alpha_1,\ldots ,\alpha_n)\in\NZ^n\}$. If there is a certain positive-degree function $d(~)$ on $\B$ (or equivalently on $A$) such that $A$ is made into an $\NZ$-graded algebra with respect to the  $\NZ$-gradation $A=\oplus_{p\in\NZ}A_p$ in which each $A_p=K\hbox{-span}\{ a^{\alpha}\in\mathcal{B}|~d(a^{\alpha})=p\}$, then we call $A$ an {\it $\NZ$-graded algebra of $(\B ,d(~))$-type}. }\v5

Let $A=K[a_1,\ldots ,a_n]$ be a finitely generated $K$-algebra with the set of generators $\{ a_1,\ldots ,a_n\}$ and the PBW basis $\mathcal{B} =\{ a^{\alpha}=a_1^{\alpha_1}\cdots
a_n^{\alpha_n}~|~\alpha =(\alpha_1,\ldots ,\alpha_n)\in\NZ^n\}$. We next explore the sufficient and necessary condition that makes $A$ into  an $\NZ$-filtered algebra of type $(\B ,d(~))$, and furthermore, we explore the associated graded structures of $A$ determined by the $\NZ$-filtration of $A$. \v5

Given an $n$-tuple $(m_1,\ldots ,m_n)\in\NZ^n$ with all $m_i>0$, condider the  positive-degree function $d(~)$ on $A$ such that
$d(a^{\alpha})=\sum_{i=1}^n\alpha_im_i$ for all $a^{\alpha}=a_1^{\alpha_1}\cdots a_n^{\alpha_n}\in\B$, in particular $d(a_i)=m_i$, $1\le i\le n$, and for a nonzero $f=\sum_{s=1}^t\lambda_sa^{\alpha (s)}\in A$ with $\lambda_s\in K^*$ and $a^{\alpha (s)}\in\B$, $d(f)=\max\{d(a^{\alpha (s)})~|~1\le s\le t\}$. Then,
with respect to $d(~)$, not only an $\NZ$-gradation may be constructed for the vector space $A$ (as we have seen from the last part), but also an $\NZ$-filtration may be constructed for the vector space $A$,  i.e., $A$ has the $\NZ$-filtration $FA=\{F_pA\}_{p\in\NZ}$, where for every $p\in\NZ$,
$$F_pA=K\hbox{-span}\{ a^{\alpha}\in\mathcal{B}|~d(a^{\alpha})\le p\}$$ which is a finite dimensional subspace of $A$, in
particular $A_0=K$, such that $F_pA\subseteq F_{p+1}A$ for all $p\in\NZ$ and $A=\cup_{p\in\NZ}F_pA$. However, in general this $\NZ$-filtration does not necessarily satisfy the condition
$$F_pAF_qA\subseteq F_{p+q}A,\quad p, q\in\NZ ,\leqno{({\rm C})}$$
namely $A$ is not necessarily an $\NZ$-filtered algebra with respect to the $\NZ$-filtration determined by the PBW basis $\B$ of $A$ and the given positive-degree function $d(~)$ on $A$ (see Example (1) given after Definition 1.4). To remedy this problem, considering the condition (C) presented above by focusing on the representation of the product $a_ja_i$ by elements in $\mathcal{B}$, where $1\le i<j\le n$, we have the following easily verified proposition. {\parindent=0pt\v5

{\bf 1.3. Proposition}  With the notation above, $A$ is a filtered algebra with respect to the $\NZ$-filtration $FA$, or equivalently, the filtration $FA$ satisfies condition (C) presented above, if and only if $A$ satisfies the condition
$$\begin{array}{l} a_ja_i\ne 0~\hbox{and}~a_ja_i=\sum^t_{s=1}\lambda_sa^{\alpha (s)}~\hbox{implies}~ d(a^{\alpha (s)})\le m_j+m_i, \\
\hbox{where}~ \lambda_s\in K^*,~a^{\alpha (s)}\in\mathcal{B} ,~1\le i<j\le n.\end{array}\leqno{(\hbox{C}')}$$\par\QED}\v5

In conclusion,  we give the formal definition of the $\NZ$-filtered algebras specified in Proposition 1.3, as follows. {\parindent=0pt\v5

{\bf  1.4. Definition}  Let $A=K[a_1,\ldots ,a_n]$ be a finitely generated $K$-algebra with the set of generators $\{ a_1,\ldots ,a_n\}$ and the PBW basis $\mathcal{B} =\{ a^{\alpha}=a_1^{\alpha_1}\cdots
a_n^{\alpha_n}~|~\alpha =(\alpha_1,\ldots ,\alpha_n)\in\NZ^n\}$. If there is a certain positive-degree function $d(~)$ on $\B$ (or equivalently on $A$) such that $A$ is made into an $\NZ$-filtered algebra with respect to the  $\NZ$-filtration $FA=\{F_pA\}_{p\in\NZ}$ where each $F_pA=K\hbox{-span}\{ a^{\alpha}\in\mathcal{B}|~d(a^{\alpha})\le p\}$, then we call $A$ an {\it $\NZ$-filtered algebra of $(\B ,d(~))$-type}. }\v5

We now give an example to illustrate Proposition 1.1 and Proposition 1.3.{\parindent=0pt\v5

{\bf Example} (1) Let $A=K[a_1,a_2,a_3]$ be the $K$-algebra generated by $\{a_1,a_2,a_3\}$ subject to the relations $$\begin{array}{l} a_1a_2=a_2a_1,\\
a_3a_1=\lambda a_1a_3-\mu a_2^2a_3-f(a_2),\\
a_3a_2=a_2a_3,\end{array}$$ where $\lambda\in K^*$, $\mu\in K$,
$f(a_2)\in K$-span$\{1,a_2,a_2^2,\ldots ,a_2^6\}$.
Then it follows from [Li4] that $A$ has the PBW basis  $\mathcal{B}' =\{
a^{\alpha}=a_1^{\alpha_1}a_2^{\alpha_2}a_3^{\alpha_3}~|~\alpha
=(\alpha_1,\alpha_2,\alpha_3)\in\NZ^3\}$.  Of course there are many choices for $f(a_2)$ and $d(~)$ such that  $A$ forms an $\NZ$-graded or $\NZ$-filtered algebra of $(\B ',d(~))$-type. For instance, consider the positive-degree function on $\B '$ such that $d(a_1)=2$, $d(a_2)=1$ and $d(a_3)=4$. Then, with $f(a_2)=a_2^6$, $A$ is an $\NZ$-graded algebra of $(\B ',d(~))$-type; noticing that with respect to the given positive-degree function $d(~)$ on $\B '$ we have $d(f(a_2))\le 6$, it turns out that $A$ is an $\NZ$-filtered algebra of $(\B ',d(~))$-type. But if we consider the positive-degree function such that $d(a_i)=1$, $1\le i\le 3$, then, $A$ is not an $\NZ$-graded algebra of $(\B ',d(~))$-type with respect to the gradation $A=\oplus_{p\in\NZ}A_p$ where each $A_p=K$-span$\{a^{\alpha}\in\B '~|~d(a^{\alpha})=p\}$, because with $\lambda\ne 0$ and $\mu\ne 0$, $a_3a_1=\lambda a_1a_3+\mu
a_2^2a_3+f(a_2)$ in which $d(a_2^2a_3)>2$, thereby $A_1A_1\not\subset A_2$; $A$ is also not an $\NZ$-filtered algebra of $(\B ',d(~))$-type with respect to the filtration $FA=\{F_pA\}_{p\in\NZ}$ where each $F_pA=K$-span$\{a^{\alpha}\in\B '~|~d(a^{\alpha})\le p\}$, because with $\lambda\ne 0$ and $\mu\ne 0$, $a_3a_1=\lambda a_1a_3+\mu a_2^2a_3+f(a_2)$ in which $d(a_2^2a_3)>2$, thereby $F_1AF_1A\not\subset F_2 A$. }\v5

Let $A$ be an $\NZ$-filtered algebra of $(\B ,d(~))$-type with the filtration $FA=\{F_pA\}_{p\in\NZ}$. Note that $FA$ is constructed with respect to a positive-degree function $d(~)$ on $A$ such that $F_pA=K$-span$\{a^{\alpha}\in\B~|~d(a^{\alpha})\le p\}$, $p\in\NZ$, in particular $F_0A=K$. It turns out that $FA$  is {\it separated} in the sense that if $f$ is a {\it nonzero} element\index{separated}
of $L$, then either $f\in F_0A=K$ or $f\in F_pA-F_{p-1}A$ for some
$p>0$. Actually, this tells us that for a nonzero $f\in A$,
$$d(f)=\left\{\begin{array}{ll} 0,&\hbox{if}~f\in F_0A=K,\\
p,&\hbox{if}~f\in F_pA-F_{p-1}A~\hbox{for
some}~p>0.\end{array}\right.$$
In order to deal with the associated graded structure of $A$ determined by $FA$ below, we first  highlight an intrinsic   property of $d(~)$ with respect to $FA$, as follows.
{\parindent=0pt\v5

{\bf 1.5. Lemma}   If $f=\sum_{s=1}^t\lambda_sa^{\alpha (s)}\in A$ with
$\lambda_s\in K^*$ and $a^{\alpha (s)}\in\B$, then $d(f)=p$ if and only if $d(a^{\alpha (s')})=p$ for some $1\le s'\le t$.\vskip 6pt

{\bf Proof} Exercise.\QED}\v5

Let $A$ be an  $\NZ$-filtered  algebra of $(\B ,d(~))$-type with filtration $FA=\{ F_pA\}_{p\in\NZ}$. Then $A$ has the {\it associated $\NZ$-graded $K$-algebra}
$G(A)=\oplus_{p\in\NZ}G(A)_p$ with $G(A)_0=F_0A=K$ and
$G(A)_p=F_pA/F_{p-1}A$ for $p\ge 1$, where for $\OV f=f+F_{p-1}A\in
G(A)_p$, $\OV g=g+F_{q-1}A$, the multiplication is given by $\OV
f\OV g=fg+F_{p+q-1}A\in G(A)_{p+q}$. Another $\NZ$-graded
$K$-algebra determined by $FA$ is the {\it Rees algebra} $\T A$ of $A$,\index{Rees algebra}
which is defined as $\T A=\oplus_{p\in\NZ}\T A_p$ with $\T
A_p=F_pA$, where the multiplication of $\T A$ is induced by
$F_pAF_qA\subseteq F_{p+q}A$, $p, q\in\NZ$. For convenience, we fix
the following notations once and for all:{\parindent=.5truecm\par

\item{$\bullet$} If $h\in G(A)_p$ and $h\ne 0$ (i.e., $h$ is a nonzero degree-$p$ homogeneous element of $G(A)$), then, as in the foregoing discussion,
we write $d_{\rm gr}(h)$ for the gr-degree of $h$ as a homogeneous
element of $G(A)$, i.e., $d_{\rm gr}(h)=p$.\par

\item{$\bullet$} If $H\in \T A_p$ and $H\ne 0$, then
we write $d_{\rm gr}(H)$ for the gr-degree of the nonzero degree-p homogeneous element $H$ of $\T A$, i.e., $d_{\rm gr}(H)=p$.\v5}

Concerning the $\NZ$-graded structure of $G(A)$, if $f\in A$ with
$d(f)=p$, then by Lemma 1.5,  the coset $f+F_{p-1}A$
represented by $f$ in $G(A)_p$ is a nonzero homogeneous element of
degree $p$. If we denote this homogeneous element by $\sigma (f)$
(in the literature it is referred to as the principal symbol of
$f$), then  $d(f)=p=d_{\rm gr}(\sigma (f))$. However,
considering the Rees algebra $\T A$ of $A$, we note that a nonzero
$f\in F_qA$ represents a homogeneous element of degree $q$ in $\T
A_q$ on one hand, and on the other hand it represents a homogeneous
element of degree $q_1$ in $\T A_{q_1}$, where $q_1=d(f)\le q$. So, for a nonzero $f\in F_pA$, we denote the
corresponding homogeneous element of degree $p$ in $\T A_p$ by
$h_p(f)$, while we use $\T f$ to denote the homogeneous element
represented by $f$ in $\T A_{p_1}$ with $p_1=d(f)\le p$.
Thus,  $d_{\rm gr}(\T f)=d(f)$, and we see that $h_p(f)=\T
f$ if and only if $d(f)=p$.
\par

Furthermore, if we write $Z$ for the homogeneous element of degree 1
in $\T A_1$ represented by the multiplicative identity element 1,
then $Z$ is a central regular element of $\T A$, i.e., $Z$ is not a
divisor of zero and is contained in the center of $\T A$. Bringing
this homogeneous element $Z$ into play,  the homogeneous elements of
$\T A$ are featured as follows:{\parindent=1truecm\vskip6pt

\item{$\bullet$} If $f\in A$ with $d(f)=p_1$ then for all $p\ge p_1$,
$h_p(f)=Z^{p-p_1}\T f$. In other words, if $H\in\T A_p$ is a nonzero
homogeneous element of degree $p$, then there is a unique $f\in F_pA$
such that $H=Z^{p-d(f)}\T f=\T f+(Z^{p-d(f)}-1)\T f$.
\par}{\parindent=0pt\vskip 6pt

With the notation above, it follows that considering the $K$-algebra homomorphism
$$\begin{array}{cccc}\VF :& \T{A}&\mapright{}{}&G(A)=\displaystyle\bigoplus_{p\in\NZ}F_pA/F_{p-1}A\\
&H&\mapsto &f+F_{p-1}A\end{array}$$
we have  $G(A)\cong\T A/\langle Z\rangle$ where $\langle Z\rangle$ is the ideal generated by $Z$ in $\T{A}$, and considering the $K$-algebra homomorphism
$$\begin{array}{cccc} \psi :&\T{A}&\mapright{}{}&A\\
&H&\mapsto &f\end{array}$$
we have $A\cong \T A/\langle 1-Z\rangle$ where $\langle 1-Z\rangle$ is the ideal generated by $1-Z$ in $\T{A}$.}\v5

The proposition  presented below will very much help us to construct graded monomial orderings (see next section for the definition) on  $\NZ$-graded and $\NZ$-filtered solvable polynomial algebras of $(\B ,d(~))$-type in the next two sections, so that the results of [Li5, Ch4, Ch5] may be better realized with such graded orderings.{\parindent=0pt\v5

{\bf 1.6. Proposition}  Let $A=K[a_1,\ldots ,a_n]$ be a finitely generated $K$-algebra with the set of generators $\{ a_1,\ldots ,a_n\}$ and the PBW basis $\mathcal{B} =\{ a^{\alpha}=a_1^{\alpha_1}\cdots
a_n^{\alpha_n}~|~\alpha =(\alpha_1,\ldots ,\alpha_n)\in\NZ^n\}$.}\par

(i) Suppose that $A$ is an $\NZ$-filtered algebra of $(\B ,d(~))$-type with the filtration $FA=\{F_pA\}_{p\in\NZ}$ where each $F_pA=K$-span$\{a^{\alpha}\in\B~|~d(a^{\alpha})\le p\}$. With the  notation as before, the following statements hold.{\parindent=1.82truecm\par

\item{(a)} For $f,g\in A$, suppose
$d(f)=p_1$ and $d(g)=p_2$.  Then $d(fg)=p_1+p_2$ if and only if $\sigma (f)\sigma (g)\ne 0$. In the case that  $\sigma (f)\sigma (g)\ne 0$ we have $\sigma (f)\sigma (g)=\sigma (fg)$ and $\T f\T g=\widetilde{fg}$. Also if $p_1+p_2\le p$, then $h_p(fg)=Z^{p-p_1-p_2}\T f\T g$.\par

\item{(b)} For $a^{\alpha}=a_1^{\alpha_1}\cdots a_n^{\alpha_n}\in\B$, we have
$$\begin{array}{l} d(a^{\alpha})=d_{\rm gr}(\sigma(a^{\alpha}))=d_{\rm gr}(\widetilde{a^{\alpha}}), \\
\sigma(a^{\alpha})=\sigma(a_1)^{\alpha_1}\cdots\sigma(a_n)^{\alpha_n}=\sigma(a)^{\alpha},~\widetilde{a^{\alpha}}=\T{a_1}^{\alpha_1}\cdots \T{a_n}^{\alpha_n}=\T a^{\alpha}. \end{array}$$

\item{(c)} $G(A)~\hbox{has the PBW basis}~\sigma (\B)=\{
\sigma (a)^{\alpha}~|~a^{\alpha}\in\B\}$ and
$\T A~\hbox{has the PBW basis}~\T{\B}=\{\T a^{\alpha}Z^m~|~a^{\alpha}\in\B,m\in\NZ\}$.\par}

(ii) Suppose that $A$ is either an $\NZ$-graded algebra of $(\B ,d(~))$-type or an $\NZ$-filtered algebra of $(\B ,d(~))$-type, and let $f=\sum_{s=1}^t\lambda_sa^{\alpha (s)}\in A$ with
$\lambda_s\in K^*$, $a^{\alpha (s)}\in\B$ and $d(f)=p$. Then by Lemma 1.5, $f$ is associated to a unique element
$$\LH_d(f)=\sum_{d(a^{\alpha(s')})=p}\lambda_{s'}a^{\alpha (s')},\quad 1\le s'\le t,$$
such that $d(f)=p=d(\LH_d(f))$. The element $\LH_d(f)$ is usually referred to as the {\it leading homogeneous part} of $f$.  In the commutative case $\LH_d(f)$ is called the {\it degree form} of $f$ (e.g. see [KR2]), and in the noncommutative case the algebra defined by leading homogeneous parts of an ideal is studied in [Li2] and [Li3]. With the notation as fixed,
the following holds true.{\parindent=1.3truecm
\item{(a)}  for any nonzero $f,g\in A$ we have $d(fg)=d(f)+d(g)$ whenever $\LH_d(f)\LH_d(g)\ne 0$; \par

\item{(b)} For $a^{\alpha},a^{\beta},a^{\eta}\in\B$, if $\LH_d(a^{\alpha}a^{\beta}a^{\eta})\not\in\{0,1\}$ and $a^{\beta}\ne
\LH_d(a^{\alpha}a^{\beta}a^{\eta})$,
then $d(a^{\beta})<d(\LH_d(a^{\alpha}a^{\beta}a^{\eta}))$;\par

\item{(c)} For $a^{\gamma},a^{\alpha},a^{\beta}, a^{\eta}\in\B$, if
$d(a^{\alpha})< d(a^{\beta})$, $\LH_d(a^{\gamma}a^{\alpha}a^{\eta})\ne
0$ and $\LH_d(a^{\gamma}a^{\beta}a^{\eta})\not\in \{ 0,1\}$, then
$d(a^{\gamma}a^{\alpha}a^{\eta})<
d(a^{\gamma}a^{\beta}a^{\eta})$.\par}

{\parindent=0pt\v5

{\bf Proof} All statements may be checked directly by referring to Proposition 1.1, Proposition 1.3  and Lemma 1.5. So we leave the details as an exercise.\QED\v5

{\bf Remark} Let $A$ be either an $\NZ$-graded algebra of $(\B ,d(~))$-type or an $\NZ$-filtered algebra of $(\B ,d(~))$-type.   Proposition 1.6(ii) may also enable us to have a graded monomial ordering on $\B$ (see next section for the definition) so that $A$ may have a Gr\"obner basis theory (though probably theoretical only, i.e., not necessarily realizable in an algorithmic way). More precisely, given a {\it total ordering} $\prec$ on $\B$ (which may not necessarily be a well-ordering),  we can define a new ordering $\prec_{gr}$ on $\B$ as follows: for $a^{\alpha},a^{\beta}\in\B$,
$$a^{\alpha}\prec_{gr}a^{\beta}\Leftrightarrow\left\{\begin{array}{l} d(a^{\alpha})<d(a^{\beta}),\\
\hbox{or}\\
d(a^{\alpha})=d(a^{\beta})~\hbox{and}~a^{\alpha}\prec a^{\beta}.\end{array}\right.$$
Clearly, the obtained $\prec_{gr}$ is now a well-ordering on $\B$ and this, in turn, determines a unique {\it leading monomial} $\LM (f)$ for every nonzero $f=\sum^t_{s=1}\lambda_sa^{\alpha (s)}\in A$, where if $d(f)=p$ and $\LH_d(f)=\sum_{d(a^{\alpha(s')})=p}\lambda_{s'}a^{\alpha (s')}$, then $\LM (f)=a^{\alpha (s')}$ for some $s'$ and thereby $d(\LM (f))=d(\LH_d(F))$. If furthermore the ordering $\prec$ satisfies the conditions (b) and (c) (with $\LH_d(~)$ replaced by the leading monomial $\LM_{\prec}(~)$), then, combining Proposition 1.6(ii), it follows that $\prec_{gr}$ is a graded monomial ordering on $\B$. } \v5

\section*{2. Solvable Polynomial algebras}

{\bf 2.1. Definition} ([K-RW], [LW]) Let $A=K[a_1,\ldots ,a_n]$ be a
finitely generated $K$-algebra. Suppose that  $A$ has the PBW
$K$-basis $\B=\{a^{\alpha}=a_1^{\alpha_1}\cdots
a_n^{\alpha_n}~|~\alpha =(\alpha_1,\ldots ,\alpha_n)\in\NZ^n\}$,
and that $\prec$ is a (two-sided) monomial ordering on $\B$. If for
all $a^{\alpha}=a_1^{\alpha_1}\cdots a_n^{\alpha_n}$,
$a^{\beta}=a_1^{\beta_1}\cdots a^{\beta_n}_n\in\B$, the following
holds:
$$\begin{array}{rcl} a^{\alpha}a^{\beta}&=&\lambda_{\alpha ,\beta}a^{\alpha +\beta}+f_{\alpha ,\beta},\\
&{~}&\hbox{where}~\lambda_{\alpha ,\beta}\in K^*,~a^{\alpha
+\beta}=a_1^{\alpha_1+\beta_1}\cdots a_n^{\alpha_n+\beta_n}~\hbox{and}~f_{\alpha ,\beta}\in K\hbox{-span}\B\\
&{~}&\hbox{such that either}~f_{\alpha ,\beta}=0~\hbox{or}~\LM
(f_{\alpha ,\beta})\prec a^{\alpha +\beta},\end{array}\leqno{(\hbox{S})}$$ where $\LM
(f_{\alpha ,\beta})$ stands for the leading monomial of $f_{\alpha
,\beta}$ with respect to $\prec$, then $A$ is called a {\it solvable
polynomial algebra}.\par

Usually $(\B ,\prec )$ is referred to an {\it admissible system} of $A$.{\parindent=0pt\v5

{\bf 2.2. Definition} Let $A$ be a solvable polynomial algebra with admissible system $(\B ,\prec )$. If $d(~)$ is a positive-degree function $d(~)$ on $\B$ such that  for all $a^{\alpha},a^{\beta}\in\B$,
$$a^{\alpha}\prec a^{\beta}~\hbox{implies}~d(a^{\alpha})<d(a^{\beta}),$$
then we call $\prec$ a {\it graded monomial ordering} with respect to $d(~)$.}\v5

The next proposition provides us with a constructive characterization of solvable polynomial algebras.{\parindent=0pt\v5

{\bf 2.3. Proposition}  [Li4, Theorem 2.1] Let $A=K[a_1,\ldots
,a_n]$ be a finitely generated $K$-algebra, and let $\KS =K\langle
X_1,\ldots ,X_n\rangle$ be the free $K$-algebras with the standard
$K$-basis $\mathbb{B}=\{ 1\}\cup\{X_{i_1}\cdots X_{i_s}~|~X_{i_j}\in
X,~s\ge 1\}$. The following two statements are equivalent:}\par

(i) $A$ is a solvable polynomial algebra in the sense of Definition
2.1.\par

(ii) $A\cong \OV A=\KS /I$ via the $K$-algebra epimorphism $\pi_1$:
$\KS \r A$ with $\pi_1(X_i)=a_i$, $1\le i\le n$, $I=$ Ker$\pi_1$,
satisfying  {\parindent=1.85truecm

\item{(a)} with respect to some monomial ordering $\prec_{_X}$ on $\mathbb{B}$, the ideal $I$ has a
finite Gr\"obner basis $G$ and the reduced Gr\"obner basis of $I$ is
of the form
$$\G =\left\{ g_{ji}=X_jX_i-\lambda_{ji}X_iX_j-F_{ji}
~\left |~\begin{array}{l} \LM (g_{ji})=X_jX_i,\\ 1\le i<j\le
n\end{array}\right. \right\} $$ where $\lambda_{ji}\in K^*$,
$\mu^{ji}_q\in K$, and
$F_{ji}=\sum_q\mu^{ji}_qX_1^{\alpha_{1q}}X_2^{\alpha_{2q}}\cdots
X_n^{\alpha_{nq}}$ with $(\alpha_{1q},\alpha_{2q},\ldots
,\alpha_{nq})\in\NZ^n$, thereby $\B =\{ \OV X_1^{\alpha_1}\OV
X_2^{\alpha_2}\cdots \OV X_n^{\alpha_n}~|~$ $\alpha_j\in\NZ\}$ forms
a PBW $K$-basis for $\OV A$, where each $\OV X_i$ denotes the coset
of $I$ represented by $X_i$ in $\OV A$; and

\item{(b)} there is a (two-sided) monomial ordering
$\prec$ on $\B$ such that $\LM (\OV{F}_{ji})\prec \OV X_i\OV X_j$
whenever $\OV F_{ji}\ne 0$, where $\OV F_{ji}=\sum_q\mu^{ji}_q\OV
X_1^{\alpha_{1i}}\OV X_2^{\alpha_{2i}}\cdots \OV X_n^{\alpha_{ni}}$,
$1\le i<j\le n$. \par}\QED\v5

\section*{3. Graded Ordering for $\NZ$-graded Solvable polynomial Algebras of $(\B ,d(~))$-type}
Let  $A=K[a_1,\ldots ,a_n]$ be a finitely generated $K$-algebra with the PBW basis $\mathcal{B} =\{ a^{\alpha}=a_1^{\alpha_1}\cdots
a_n^{\alpha_n}~|~\alpha =(\alpha_1,\ldots ,\alpha_n)\in\NZ^n\}$. Recall from Definition 1.2 that $A$ is called an $\NZ$-graded algebra of $(\B ,d(~))$-type if there is a certain positive-degree function $d(~)$ on $\B$ (or equivalently on $A$) such that $A$ is made into an $\NZ$-graded algebra with respect to the  $\NZ$-gradation $A=\oplus_{p\in\NZ}A_p$ in which each $A_p=K\hbox{-span}\{ a^{\alpha}\in\mathcal{B}~|~d(a^{\alpha})=p\}$. By Proposition 1.1(i), this is also equivalent to say that $A=\oplus_{p\in\NZ}A_p$ in which every degree-$p$ homogeneous part $A_p$ is a subspace of $A$ such that $A_pA_q\subset A_{p+q}$ for all $p,q\in\NZ$, and moreover, \par

(a) $A_0=K$, and\par

(b) if $a^{\alpha}\in\mathcal{B}$, then $a^{\alpha}\in A_q$ for some $q\in\NZ$.\par{\parindent=0pt

We recall also from the definition of a solvable polynomial algebra (Definition 2.1) that the condition (S) is equivalent to
$$\begin{array}{rcl} a_ja_i&=&\lambda_{ji}a_ia_j+f_{ji}\\
&{~}&\hbox{where}~\lambda_{ji}\in K^*, ~f_{ji}=\sum\mu_ka^{\alpha
(k)}\in K\hbox{-span}\B\\
&{~}&\hbox{with}~\LM (f_{ji})\prec a_ia_j~\hbox{if}~f_{ji}\ne
0,~1\le j<i\le n.\end{array}\leqno{(\hbox{S}')}$$

Now,  it follows from Proposition 1.1(ii) that we first have the following\v5

{\bf 3.1. Proposition}  Let $A=K[a_1,\ldots, a_n]$ be a solvable polynomial algebra with admissible system $(\B ,\prec )$. The following statements are equivalent.}\par
(i)  $A$ is an $\NZ$-graded algebra of $(\B ,d(~))$-type (in the sense of Definition 1.2).\par

(ii) There is a  positive-degree function $d(~)$ on $A$  such that for $1\le i<j\le n$, all the relations
$a_ja_i=\lambda_{ji}a_ia_j+f_{ji}$ with $f_{ji}=\sum\mu_ka^{\alpha
(k)}$ presented in (S$'$) above satisfy  $d(a^{\alpha
(k)})=d(a_ia_j)=d(a_j)+d(a_i)$ whenever $\mu_k\ne 0$.\par\QED\v5

The proposition (or more precisely its proof) below shows us that every $\NZ$-graded solvable polynomial algebra of $(\B ,d(~))$-type has a graded monomial ordering $\prec_{gr}$ with respect to the same positive-degree function $d(~)$ on $\B$ as described in  the remark made at the end of Section 1.{\parindent=0pt\v5

{\bf 3.2. Proposition}  Let $A=K[a_1,\ldots ,a_n]$ be a finitely generated $K$-algebra with the PBW basis $\mathcal{B} =\{ a^{\alpha}=a_1^{\alpha_1}\cdots
a_n^{\alpha_n}~|~\alpha =(\alpha_1,\ldots ,\alpha_n)\in\NZ^n\}$.  The following statements are equivalent.}\par

(i) $A$ is a solvable polynomial algebra with a graded monomial ordering $\prec_{gr}$ with respect to a certain positive-degree function $d(~)$ on $\B$  such that for $1\le i<j\le n$, all the relations
$a_ja_i=\lambda_{ji}a_ia_j+f_{ji}$ with $f_{ji}=\sum\mu_ka^{\alpha
(k)}$ presented in (S2$'$) above satisfy  $d(a^{\alpha
(k)})=d(a_ia_j)=d(a_j)+d(a_i)$ whenever $\mu_k\ne 0$;\par

(ii) $A$ is an $\NZ$-graded solvable polynomial algebra of $(\B ,d(~))$-type. {\parindent=0pt\vskip 6pt

{\bf Proof} (i) $\Rightarrow$ (ii) This follows immediately from Proposition 3.1.\par

(ii) $\Rightarrow$ (i) Let $A$ be an $\NZ$-graded solvable polynomial algebra of $(\B ,d(~))$-type. By Proposition 3.1 it remains to prove that $A$ has a graded monomial ordering with respect to the given positive-degree function $d(~)$. As $A$ is a solvable polynomial algebra, there is a monomial ordering $\prec$ on $\B$. If we define a new ordering $\prec_{gr}$ on $\B$ subject to the rule: for $a^{\alpha}$, $a^{\beta}\in\B$,
$$a^{\alpha}\prec_{gr}a^{\beta}\Leftrightarrow\left\{\begin{array}{l} d(a^{\alpha})<d(a^{\beta})\\
\hbox{or}\\
d(a^{\alpha})=d(a^{\beta})~\hbox{and}~a^{\alpha}\prec a^{\beta},\end{array}\right.$$
then, noticing that $A$ is a domain and using Proposition 1.6(ii), it is straightforward to check that $\prec_{gr}$ is a  graded monomial ordering on $\B$ with respect to $d(~)$, as desired.\QED}\v5

As we will see in the next section, that every  $\NZ$-filtered solvable polynomial algebra of $(\B ,d(~))$-type has a graded monomial ordering $\prec_{gr}$. This enables us to show that the associated  graded algebra $G(A)$ of an $\NZ$-filtered solvable
polynomial algebra $A$ of $(\B ,d(~))$-type  is an $\NZ$-graded solvable polynomial algebra  of $(\sigma(\B ) ,d(~))$-type (thereby it has a graded monomial ordering), and the Rees algebra $\T{A}$ of $A$ is an $\NZ$-graded solvable polynomial algebra of $(\T{\B},d(~))$-type (thereby it has a graded monomial ordering).\v5

Finally, let us also have a review of the algebra $A$ given in Example (1) of Section 1. It follows from [Li4, Proposition 3.1 and Proposition 3.2]  that with $f(a_2)=a_2^6$ and the positive-degree function such that $d(a_1)=2$, $d(a_2)=1$ and $d(a_3)=4$ , $A$ is turned into an $\NZ$-graded solvable polynomial algebra of $(\B ',d(~))$-type with the graded monomial ordering $a_3\prec_{grlex}a_2\prec_{grlex}a_1$.\v5

\section*{4. Graded Ordering for $\NZ$-filtered Solvable polynomial Algebras of $(\B ,d(~))$-type}

Let  $A=K[a_1,\ldots ,a_n]$ be a finitely generated $K$-algebra with the PBW basis $\mathcal{B} =\{ a^{\alpha}=a_1^{\alpha_1}\cdots
a_n^{\alpha_n}~|~\alpha =(\alpha_1,\ldots ,\alpha_n)\in\NZ^n\}$. Recall from Definition 1.4 that $A$ is called an $\NZ$-filtered algebra of $(\B ,d(~))$-type if there is a certain positive-degree function $d(~)$ on $\B$ (or equivalently on $A$) such that $A$ is made into an $\NZ$-filtered  algebra with respect to the  $\NZ$-filtration $FA=\{F_pA\}_{p\in\NZ}$ in which each $F_pA=K\hbox{-span}\{ a^{\alpha}\in\mathcal{B}~|~d(a^{\alpha})\le p\}$. By Proposition 1.3, this is also equivalent to say that $A$ satisfies the condition\par
\centerline{$(\hbox{C}')\quad\quad\begin{array}{l} a_ja_i\ne 0~\hbox{and}~a_ja_i=\sum^t_{s=1}\lambda_sa^{\alpha (s)}~\hbox{implies}~ d(a^{\alpha (s)})\le m_j+m_i, \\
\hbox{where}~ \lambda_s\in K^*,~a^{\alpha (s)}\in\mathcal{B} ,~1\le i<j\le n.\end{array}$}
{\parindent=0pt\par
 Also recall that in the last section we have pointed out that the condition (S) of Definition 2.1 is equivalent to the condition
$$\begin{array}{rcl} a_ja_i&=&\lambda_{ji}a_ia_j+f_{ji}\\
&{~}&\hbox{where}~\lambda_{ji}\in K^*, ~f_{ji}=\sum\mu_ka^{\alpha
(k)}\in K\hbox{-span}\B\\
&{~}&\hbox{with}~\LM (f_{ji})\prec a_ia_j~\hbox{if}~f_{ji}\ne
0,~1\le j<i\le n.\end{array}\leqno{(\hbox{S}')}$$
Now,  it follows from Proposition 1.3 that we first have the following{\parindent=0pt\v5

{\bf 4.1. Proposition}  Let $A=K[a_1,\ldots, a_n]$ be a solvable polynomial algebra with admissible system $(\B ,\prec )$. The following statements are equivalent.}}\par
(i)  $A$ is an $\NZ$-filtered algebra of $(\B ,d(~))$-type (in the sense of Definition 1.4).\par

(ii) There is a  positive-degree function $d(~)$ on $A$  such that for $1\le i<j\le n$, all the relations
$a_ja_i=\lambda_{ji}a_ia_j+f_{ji}$ with $f_{ji}=\sum\mu_ka^{\alpha
(k)}$ presented in (S$'$)  satisfy  $d(a^{\alpha
(k)})\le d(a_ia_j)=d(a_j)+d(a_i)$ whenever $\mu_k\ne 0$.\par\QED\v5

The proposition (or more precisely its proof) below shows us that  every $\NZ$-filtered solvable polynomial algebra of $(\B ,d(~))$-type has a graded monomial ordering $\prec_{gr}$ with respect to the same positive-degree function $d(~)$ on $\B$ as described in  the remark made at the end of Section 1.{\parindent=0pt\v5

{\bf 4.2. Proposition}   Let $A=K[a_1,\ldots ,a_n]$ be a finitely generated $K$-algebra with the PBW basis $\mathcal{B} =\{ a^{\alpha}=a_1^{\alpha_1}\cdots
a_n^{\alpha_n}~|~\alpha =(\alpha_1,\ldots ,\alpha_n)\in\NZ^n\}$.  The following statements are equivalent.}\par

(i) $A$ is a solvable polynomial algebra with a graded monomial ordering $\prec_{gr}$ with respect to a certain positive-degree function $d(~)$ on $\B$; \par

(ii) $A$ is an $\NZ$-filtered solvable polynomial algebra of $(\B ,d(~))$-type. {\parindent=0pt\vskip 6pt

{\bf Proof} (i) $\Rightarrow$ (ii) If $\prec_{gr}$ is a graded monomial ordering with respect to the positive-degree function $d(~)$ on $\B$, then by  the condition (S$'$) mentioned above we have for $1\le i<j\le n$, all the relations
$a_ja_i=\lambda_{ji}a_ia_j+f_{ji}$ with $f_{ji}=\sum\mu_ka^{\alpha
(k)}$  satisfy  $d(a^{\alpha
(k)})\le d(a_ia_j)=d(a_j)+d(a_i)$ whenever $\mu_k\ne 0$. Thus it follows from Proposition 4.1 that $A$ is an $\NZ$-filtered solvable polynomial algebra of $(\B ,d(~))$-type.}\par

(ii) $\Rightarrow$ (i) Let $A$ be an $\NZ$-filtered solvable polynomial algebra of $(\B ,d(~))$-type, where $d(~)$ is a  positive-degree function on $\B$. In order to reach (i),  we proceed  to show that $A$ has a graded monomial ordering with respect to the given positive-degree function $d(~)$ on $\B$. To this end, note that $A$ is a solvable polynomial algebra, thereby there is a monomial ordering $\prec$ on $\B$. If we define a new ordering $\prec_{gr}$ on $\B$ subject to the rule: for $a^{\alpha}$, $a^{\beta}\in\B$,
$$a^{\alpha}\prec_{gr}a^{\beta}\Leftrightarrow\left\{\begin{array}{l} d(a^{\alpha})<d(a^{\beta})\\
\hbox{or}\\
d(a^{\alpha})=d(a^{\beta})~\hbox{and}~a^{\alpha}\prec a^{\beta},\end{array}\right.$$
then, noticing that $A$ is a domain and using Proposition 1.6(ii), it is straightforward to check that $\prec_{gr}$ is a  graded monomial ordering on $\B$ with respect to $d(~)$, as desired.\QED\v5

With the aid of Proposition 4.1, the following examples may be better understood.{\parindent=0pt\v5

{\bf Example} (1) If $A=K[a_1,\ldots ,a_n]$ is an $\NZ$-graded
solvable polynomial algebra of $(\B ,d(~))$-type, i.e., $A=\oplus_{p\in\NZ}A_p$ with the degree-$p$
homogeneous part $A_p=K$-span$\{a^{\alpha}\in\B~|~d(a^{\alpha})=p\}$, then,  $A$ is turned into an
$\NZ$-filtered solvable polynomial algebra of $(\B ,d(~))$-type by the same
positive-degree function $d(~)$ and the $\NZ$-filtration
$FA=\{F_pA\}_{p\in\NZ}$ with each $F_pA=\oplus_{q\le p}A_q$.} \v5

The next example  provides $\NZ$-filtered solvable polynomial
algebras of $(\B ,d(~))$-type in which some generators may have degree $\ge
2$.{\parindent=0pt\v5

{\bf Example} (2)  Consider Example (1) given in Section 1, in which the algebra $A=K[a_1,a_2,a_3]$ has the PBW basis  $\mathcal{B}' =\{ a^{\alpha}=a_1^{\alpha_1}a_2^{\alpha_2}a_3^{\alpha_3}~|~\alpha
=(\alpha_1,\alpha_2,\alpha_3)\in\NZ^3\}$, $a_1a_2=a_2a_1$, $a_3a_2=a_2a_3$, and $a_3a_1=\lambda a_1a_3+\mu a_2^2a_3+f(a_2)$ with $f(a_2)\in K$-span$\{ 1,a_2,a_2^2,\ldots ,a_2^6\}$. Then by  Proposition 3.1, we see that with  the positive-degree function such that $d(a_1)=2$, $d(a_2)=1$ and $d(a_3)=4$ , $A$ is turned into an $\NZ$-filtered  solvable polynomial algebra of $(\B ',d(~))$-type.}\v5

Let $A=K[a_1,\ldots ,a_n]$  be an  $\NZ$-filtered solvable polynomial algebra of $(\B ,d(~))$-type with the filtration $FA=\{F_pA\}_{p\in\NZ}$ in which each $F_pA=K$-span$\{ a^{\alpha}\in\B~|~d(a^{\alpha})\le p\}$, and let $\prec_{gr}$ be a graded monomial ordering on $\B$ with respect to the same positive-degree function $d(~)$ on $\B$ (existence guaranteed  by Proposition 4.2). Considering the associated graded algebra $G(A)=\oplus_{p\in\NZ}F_pA/F_{p-1}$ of $A$ as well as the Rees algebra $\T A=\oplus_{p\in\NZ}F_pA$ of $A$ (see Section 1), we  are ready to present  the following{\parindent=0pt\v5

{\bf 4.3. Theorem}  With the notation as in Section 1 and the notation fixed above, the following statements hold.}\par

(i) The $\NZ$-graded $K$-algebra $G(A)$ is generated by $\{\sigma (a_1),\ldots ,\sigma (a_n)\}$, i.e., $G(A)=K[\sigma (a_1),\ldots ,\sigma (a_n)]$,  and $G(A)$ has the PBW basis
$$\sigma (\B)=\{
\sigma (a)^{\alpha}=\sigma (a_1)^{\alpha_1}\cdots \sigma
(a_n)^{\alpha_n}~|~\alpha =(\alpha_1,\ldots ,\alpha_n)\in\NZ^n\} .$$
By referring to the relation $a^{\alpha}a^{\beta}=\lambda_{\alpha ,\beta}a^{\alpha
+\beta}+f_{\alpha ,\beta}$ in Definition 2.1, where $\lambda_{\alpha ,\beta}\in K^*$, for  $\sigma
(a)^{\alpha}$, $\sigma (a)^{\beta}\in\sigma (\B )$, if $f_{\alpha ,\beta}=0$ then
$$\sigma (a)^{\alpha}\sigma(a)^{\beta}=
\lambda_{\alpha ,\beta}\sigma (a)^{\alpha
+\beta},~\hbox{where}~\sigma (a)^{\alpha +\beta}=\sigma
(a_1)^{\alpha_1+\beta_1}\cdots \sigma (a_n)^{\alpha_n+\beta_n};$$
in the case where $f_{\alpha ,\beta}=\sum_j\mu^{\alpha
,\beta}_ja^{\alpha (j)}\ne 0$ with $\mu^{\alpha ,\beta}_j\in K$,
$$\sigma (a)^{\alpha}\sigma(a)^{\beta}= \lambda_{\alpha ,\beta}\sigma
(a)^{\alpha +\beta}+\displaystyle{\sum_{d(a^{\alpha
(k)})=d(a^{\alpha +\beta})}}\mu^{\alpha ,\beta}_j\sigma (a)^{\alpha
(k)}.$$ Moreover, the ordering $\prec_{_{G(A)}}$ defined on $\sigma
(\B )$ subject to the rule: for $\sigma (a)^{\alpha},\sigma (a)^{\beta}\in\sigma (\B )$,
$$\sigma (a)^{\alpha}\prec_{_{G(A)}} \sigma (a)^{\beta}\Leftrightarrow a^{\alpha}\prec_{gr}a^{\beta},$$  is a graded monomial ordering
with respect to the positive-degree function $d(~)$ on $\sigma (\B )$ defined by assigning   $d(\sigma (a_i))=d(a_i)$  for $1\le i\le n$. Thus, with the data $(\sigma (\B ),\prec_{_{G(A)}}, d(~))$,  $G(A)$ is turned
into an $\NZ$-graded  solvable polynomial algebra of $(\sigma (\B ),d(~))$-type.\par

(ii) The $\NZ$-graded $K$-algebra $\T{A}$ is generated by $\{\T a_1,\ldots ,\T a_n, Z\}$, i.e., $\T A=K[\T a_1,\ldots ,\T a_n, Z]$ where $Z$ is the central
regular element of degree 1 in $\T A_1$ represented by 1,  and $\T A$
has the PBW basis $$\T{\B}=\{\T a^{\alpha}Z^m=\T
a_1^{\alpha_1}\cdots \T a_n^{\alpha_n}Z^m~|~\alpha =(\alpha_1,\ldots
,\alpha_n )\in\NZ^n,m\in\NZ\} .$$
By referring to the relation
$a^{\alpha}a^{\beta}=\lambda_{\alpha ,\beta}a^{\alpha
+\beta}+f_{\alpha ,\beta}$ in Definition 2.1, where $\lambda_{\alpha ,\beta}\in K^*$, for  $\T a^{\alpha}Z^s$, $\T a^{\beta}Z^t\in\T{\B}$,  if $f_{\alpha ,\beta}=0$ then
$$\T a^{\alpha}Z^s\cdot\T a^{\beta}Z^t=
\lambda_{\alpha ,\beta}\T a^{\alpha +\beta}Z^{s+t},~\hbox{where}~\T
a^{\alpha +\beta}=\T a_1^{\alpha_1+\beta_1}\cdots \T
a_n^{\alpha_n+\beta_n};$$ in the case where $f_{\alpha
,\beta}=\sum_j\mu^{\alpha ,\beta}_ja^{\alpha (j)}\ne 0$ with
$\mu^{\alpha ,\beta}_j\in K$,
$$\begin{array}{rcl} \T a^{\alpha}Z^s\cdot\T a^{\beta}Z^t&=&
\lambda_{\alpha ,\beta}\T a^{\alpha +\beta}Z^{s+t}+\displaystyle\sum_j\mu^{\alpha
,\beta}_j\T a^{\alpha (j)}Z^{q-m_j},\\
&{~}&\hbox{where}~q=d(a^{\alpha +\beta})+s+t,~m_j=d(a^{\alpha
(j)}).\end{array}$$ Moreover, the ordering $\prec_{_{\T A}}$ defined
on $\T{\B}$ subject to the rule: for $\T a^{\alpha}Z^s , \T a^{\beta}Z^t\in\T{\B}$,
$$\T a^{\alpha}Z^s\prec_{_{\T A}}\T a^{\beta}Z^t\Leftrightarrow\left\{\begin{array}{l}  a^{\alpha}\prec_{gr}a^{\beta}\\
\hbox{or}\\
a^{\alpha}=a^{\beta}~\hbox{and}~s<t,\end{array}\right.$$ is a monomial ordering on $\T{\B}$
(which is not necessarily a graded monomial ordering). Thus, with the data $(\T{\B},\prec_{_{\T A}}, d(~))$, where $d(~)$ is the positive-degree function on $\T A$ defined by assigning $d(Z)=1$
and $d(\T{a_i})=d (a_i)$ for $1\le i\le n$,
$\T A$ is turned into an $\NZ$-graded solvable polynomial algebra of $(\T{\B} ,d(~))$-type. {\parindent=0pt\vskip 6pt

{\bf Proof} The results stated in the theorem are just analogues of those results
concerning quadric solvable polynomial algebras established in ([LW, Section 3],
[Li1, CH.IV]), so they may be derived in a similar way as in loc. cit., or else  they may also be proved directly  by use Proposition 1.6 and thus we leave the detailed argumentation as an exercise.\QED }\v5

Let $A$ be an $\NZ$-filtered algebra of $(\B ,d(~))$-type with the filtration $FA=\{F_pA\}_{p\in\NZ}$ in which each $F_pA=K$-span$\{a^{\alpha}\in\B~|~d(a^{\alpha})\le p\}$.
We end this section by a lemma that  may be very helpful in understanding [Li6, Ch5].  {\parindent=0pt \v5

{\bf 4.4. Lemma}   With the notation as in
Theorem 4.3 and  Section 1,  if $f=\lambda a^{\alpha}+\sum_j\mu_ja^{\alpha (j)}$ with $d(f)=p$ and $\LM (f)=a^{\alpha}$, then {\parindent=1.24truecm\vskip 6pt

\item{(a)} $p=d(f)=d_{\rm gr}(\sigma (f))=d_{\rm gr}(\T f)$;\vskip 6pt

\item{(b)} $\sigma (f)=\lambda\sigma (a)^{\alpha}+\displaystyle\sum_{d(a^{\alpha (j_k)})=p}\mu_{j_k}\sigma (a)^{\alpha (j_k)}$; \par
$\LM(\sigma (f))=\sigma (a)^{\alpha}=\sigma (\LM (f))$;\vskip 6pt

\item{(c)} $\T f=\lambda\T a^{\alpha}+\displaystyle\sum_j\mu_j\T a^{\alpha (j)}Z^{p-d(a^{\alpha (j)})}$;\par
$\LM (\T f)=\T a^{\alpha}=\widetilde{\LM (f)},$\vskip 6pt}{\parindent=0pt

where $\LM (f)$, $\LM (\sigma (f))$ and $\LM (\T f)$ are taken with
respect to $\prec_{gr}$, $\prec_{_{G(A)}}$ and $\prec_{_{\T A}}$
respectively.}\vskip 6pt

{\bf Proof} All properties may be verified by referring to Section 1, in particular,  Proposition 1.6. We leave the detailed argumentation  as an exercise.\QED}
\v5
\centerline{\bf References} {\parindent=.8truecm\vskip 8pt

\item{[K-RW]} A. Kandri-Rody, V.~Weispfenning, Non-commutative
Gr\"obner bases in algebras of solvable type.  {\it J. Symbolic
Comput.}, 9(1990), 1--26.\par

\item{[Li1]} H. Li, {\it Noncommutative Gr\"obner Bases and
Filtered-Graded Transfer}. Lecture Notes in Mathematics, Vol. 1795,
Springer, 2002.\par

\item{[Li2]} H. Li, $\Gamma$-leading homogeneous algebras and Gr\"obner
bases. In: {\it Recent Developments in Algebra and Related Areas}
(F. Li and C. Dong eds.), Advanced Lectures in Mathematics, Vol. 8,
International Press \& Higher Education Press, Boston-Beijing, 2009,
155--200.
\item{[Li3]} H. Li, {\it Gr\"obner Bases in Ring Theory}. World Scientific Publishing Company, 2011.\par

\item{[Li4]} H. Li, A Note on solvable polynomial algebras. {\it Computer Science Journal of Moldova}, vol.22, 1(64), 2014, 99 -- 109.
Available at  arXiv:1212.5988[math.RA].

\item{[Li5]} H. Li, On Computation of minimal free resolutions over solvable
polynomial algebras. {\it Commentationes
Mathematicae Universitatis Carolinae}, 4(56)(2015), 447--503.

\item{[Li6]} H. Li, Notes on \GR bases and free resolutions of modules over solvable polynomial algebras).
arXiv:1510.04381v1[math.RA], 131 pages.

\item{[LW]} H. Li.  Y. Wu, Filtered-graded transfer of
Gr\"obner basis computation in solvable polynomial algebras. {\it
Communications in Algebra}, 28(1), 2000, 15--32.\par

\end{document}